\documentclass[11pt,bezier]{article}
\usepackage{amsmath,graphicx,amssymb,amsfonts,color}

\font\bg=cmbx10 scaled\magstep1 
\font\small=cmr8

\newtheorem{newlemma}{{\bf Lemma}}

\newtheorem{newteorem}{{\bf Theorem}}

\newtheorem{newkorolari}{{\bf Corollary}}

\newtheorem{newkonjek}{{\bf Conjecture}}

\newtheorem{newdefine}{{\bf Definition}}

\begin{document}
\tolerance=10000 \baselineskip18truept
\newbox\thebox
\global\setbox\thebox=\vbox to 0.2truecm{\hsize
0.15truecm\noindent\hfill}
\def\boxit#1{\vbox{\hrule\hbox{\vrule\kern0pt
     \vbox{\kern0pt#1\kern0pt}\kern0pt\vrule}\hrule}}
\def\qed{\lower0.1cm\hbox{\noindent \boxit{\copy\thebox}}\bigskip}
\def\ss{\smallskip}
\def\ms{\medskip}
\def\nt{\noindent}

\centerline{\Large\bf Some new results on domination roots of a
graph}
 \vspace{.3cm}

\bigskip

\baselineskip12truept \centerline{\bg Saeid Alikhani$^{a,b}$
}\baselineskip20truept \centerline{$^a$Department of Mathematics,
 Yazd University} \vskip-8truept \centerline{ 89195-741, Yazd,
Iran} \centerline{ $^b$School of Mathematics, Institute for
Research in Fundamental Sciences (IPM)} \centerline{P.O. Box:
19395-5746, Tehran,
             Iran.}

 \centerline{\texttt{E-mail:alikhani@yazduni.ac.ir}}
 \nt\rule{16cm}{0.1mm}

\nt{\bg ABSTRACT}

\nt Let $G$ be a simple graph of order $n$.  The domination
polynomial of $G$ is the polynomial $D(G,\lambda)=\sum_{i=0}^{n}
d(G,i) \lambda^{i}$, where $d(G,i)$ is the number of dominating
sets  of $G$ of size $i$. Every root of $D(G,\lambda)$ is called
 the domination root of $G$. We present families of graphs whose their domination polynomial have no nonzero real
      roots. We observe that these graphs have complex domination roots with positive
real part. Then, we consider the lexicographic product of two
graphs and obtain a formula for
  domination polynomial of this product.  Using this product, we construct a family of graphs which
   their domination roots are dense in all of $\mathbb{C}$.

\nt{\bf  Mathematics Subject Classification:} {\small 05C60.}

\nt\rule{16cm}{0.1mm}

\end{document}